\documentclass[12pt]{article}
\usepackage{amssymb}
\usepackage{amscd}
\usepackage{amsmath}
\usepackage{amstext}
\usepackage[all]{xy}

\title{\Large \bf Canonical Filtrations of Gorenstein Injective Modules
\thanks{2000 { Mathematics Subject Classification}. 13D07, 16E30 .}
}
\author{Edgar E. Enochs$^{1}$ and Zhaoyong Huang$^{2}$\\
{\footnotesize $^{1}$Department of Mathematics, University of Kentucky,
Lexington, KY 40506, USA;}\\ {\footnotesize $^{2}$Department of Mathematics,
Nanjing University, Nanjing 210093, Jiangsu Province P.R. China}
\\{\footnotesize (email: enochs@ms.uky.edu; huangzy@nju.edu.cn)}}
\date{}
\begin{document}
\baselineskip=18pt \maketitle
\begin{abstract} The principle ``Every result in classical homological algebra
should have a counterpart in Gorenstein homological algebra" is given in [3].
There is a remarkable body of evidence supporting this claim (cf. [2] and
[3]). Perhaps one of the most glaring exceptions is provided by the fact
that tensor products of Gorenstein projective modules need not be Gorenstein
projective, even over Gorenstein rings. So perhaps it is surprising that
tensor products of Gorenstein injective modules over Gorenstein rings
of finite Krull dimension are Gorenstein injective.\\
Our main result  is in support of the principle. Over commutative, noetherian rings injective modules have direct sum decompositions into indecomposable
modules. We will show that Gorenstein injective modules over
 Gorenstein rings of finite Krull dimension
 have filtrations analogous to those provided by these decompositions.
This result will then provide us with the tools to prove that all tensor
 products of Gorenstein injective modules over these rings are Gorenstein injective.

\end{abstract}

\vspace{0.5cm}

{\bf 1. Introduction}

\vspace{0.1cm}

Throughout this paper $R$ will denote a commutative and noetherian ring and
$Spec(R)$ will denote the set of its prime ideals.
The term module will then mean an $R$-module. An injective envelope of the
module $M$ will be denoted by $E(M)$ and
\begin{center}
 {$ 0\rightarrow M\rightarrow E^0(M)\rightarrow \cdots \rightarrow
E^n(M)\rightarrow \cdots $}
\end{center}
will denote a minimal injective resolution of $M$.

\vspace{0.2cm}

We will now give several definitions and results. For ease in appealing to
these later in the paper they will be numbered.
\begin{itemize}
\item[(1)] Every injective module is uniquely up to isomorphism the direct
sum of modules each of which is isomorphic to $E(R/P)$ for some $P\in Spec(R)$
([7], Theorem 2.5 and Proposition 3.1).
\item[(2)] We say $R$ is a Gorenstein ring if $inj.dim_{R_P}R_P < \infty $
for each $P\in Spec(R)$. If in fact $inj.dim_R R<\infty $ then $R$ is
Gorenstein and the Krull dimension of $R$ equals $inj.dim_R R$ ([1],
Corollary 3.4).
\item[(3)] If $inj.dim_R R<\infty $ (and so $R$ is Gorenstein) and if
$0\rightarrow R\rightarrow E^0(R)\rightarrow \cdots\rightarrow E ^n(R)\rightarrow 0$
is a minimal injective resolution of $R$ as a module, then for
$0\leq k \leq n$, $E ^k(R)=\oplus_{ht(P)=k} E(R/P)$ where these $P$
are in $Spec(R)$ ([1], Proposition 3.6). Furthermore when $ht(P)=k$ we have
$flat.dim E(R/P)=k$ ([8], Theorem 5.1.2).
\item[(4)] If $R$ is Gorenstein and $E,E'$ are injective modules, then
for any $k\geq 0$ the module $Tor_k(E,E')$ is  injective. \\More precisely
if $P,Q\in Spec(R)$ then $Tor_k(E(R/P), E(R/Q))=0$ unless both $P=Q$
and $k=ht(P)$. And in this case we have $Tor_k(E(R/P), E(R/P))\cong E(R/P)$
([5], Lemma 2.1 and Theorem 4.1). So using (1) we see that\\
 $Tor_k( E(R/P), E)=0$ when $E$ is injective and $k\neq ht(P)$.
\item[(5)] If $P\in Spec(R)$ a module $S$ will be said to have property
$t(P)$ if for each $r\in R-P$ we have $S\stackrel{r}{\rightarrow} S$ is
an isomorphism and if for each $x\in S$ we have $P^mx=0$ for some $m\geq 1$.
If $S$ has property $t(P)$ and property $t(Q)$ with $P\neq Q$ then it
is easy to see that $S=0$.
If $S$ has property $t(P)$ and if $N$ is any module, then $Tor_k(S,N)$ has
property $t(P)$ for all $k\geq 0$. This can be seen by using a projective
resolution of $N $ to compute this $Tor$. Consequently, if $S$ has property
$t(P)$ and $T$ has property $t(Q)$ where $P\neq Q$ we get $Tor_k(S,T)=0$
for all $k\geq 0$. For any $P\in Spec(R)$ the module $E(R/P)$ has property
$t(P)$ ([7], Lemma 3.2).
\item[(6)] We now argue that if $S$ has property $t(P)$, then so does $E(S)$.
By (1) above $E(S)$ is a direct sum of copies of $E(R/Q)$ for various
$Q\in Spec(R)$.
If $r\in R-P$ then since $S\stackrel{r}{\rightarrow}S $ is an isomorphism,
so is $E(S)\stackrel{r}{\rightarrow} E(S)$. Now
assume that $E(R/Q)$ is a summand of $E(S)$. Then for $r\in R-P$ we have
$E(R/Q)\stackrel{r}{\rightarrow} E(R/Q)$ is an isomorphism. Hence $r\in R-Q$.
So we get $Q\subset P$. We want to argue that $Q=P$. If not, let $r\in P-Q$. The extension $S\rightarrow E(S)$ is essential,
so the module $S'=E(R/Q)\cap S$ is non-zero. Let $x\in S'$, $x\neq 0$. Then
since $x\in S$ and since $S$ has property $t(P)$ we get $P^m x=0$ for
some $m\geq 1$. So $r^m x=0$. But $E(R/Q)$ has property $t(Q)$ and $r\in R-Q$.
Hence $E(R/Q)\stackrel{r}{\rightarrow} E(R/Q)$ is an isomorphism. But
then since $S'\subset E(R/Q)$ we get $S'\stackrel{r}{\rightarrow} S'$ is
an injection. But this is not possible if $r^mx=0$ where $x\in S'$ and
$x\neq 0$. So we get $Q=P$.\\
So $E(S)$ is a direct sum of copies of $E(R/P)$ and so by (5) we see
that $E(S)$ has property $t(P)$. But then the quotient module $E(S)/S$ will
have property $t(P)$. So continuing we see that all the terms $E^i(S)$ ($i\geq
1$) in a minimal injective resolution of $S$ have property $t(P)$.
\item[(7)] If $S$ has property $t(P)$ and $T$ has property $t(Q)$ and
if $P\not\subset Q$, then $Hom(S,T)=0$. For if $r\in P-Q $ and if
$f(x)=y$ for some $f\in Hom(S,T)$ we have $r^nx=0$ for some $n\geq 1$ and
so $r^ny=0$. But since $r\not\in Q$ this is possible only if $y=0$. So we
get $f=0$.
\item[(8)] We recall that a module $G$ is said to be Gorenstein injective
if and only if there is an exact sequence
\begin{center}
{$ \cdots \rightarrow E_2\rightarrow E_1 \rightarrow E_0\rightarrow
E^0\rightarrow E^1\rightarrow E^2\rightarrow \cdots $}
\end{center}
of injective modules with $G=Ker(E^0\rightarrow E^1) $ and such that
$Hom(E,-)$ leaves the sequence exact whenever $E$ is an injective module.
For the rest of (8) we assume that
 $R$ is a Gorenstein ring of finite Krull dimension $n$. If
$n\geq 1$, a module $G$ is Gorenstein injective if and only if there
is an exact sequence
\begin{center}
{$E_{n-1}\rightarrow \cdots \rightarrow E_1\rightarrow E_0\rightarrow G\rightarrow 0$}
\end{center}
with $E_{n-1},\cdots, E_0$ injective modules. This result gives that
the class of Gorenstein injective modules over such $R$ is closed under arbitrary direct sums. Also if $n=0$ then every module
$G$ is Gorenstein injective (see [4], Theorem 4.2 for both these claims). As
a consequence we get that if $P$ is a minimal prime ideal of $R$
and if $G$ is an $R_P$-module, then $G$ is a Gorenstein injective $R$-module.
This follows from the observation that $R_P$ is a flat $R$-module, so any
injective $R_P$-module is also an injective $R$-module. Hence an exact
sequence $\cdots\rightarrow E_2\rightarrow E_1\rightarrow E_0\rightarrow G\rightarrow 0$
of $R_P$-modules with the $E_k$ injective $R_P$-modules gives us an
exact sequence of $R$-modules with the $E_k$ injective $R$-modules. \\We need a slightly
stronger version of the result above. So again we suppose $R$ is Gorenstein
and of Krull dimension $n$ but with $n\geq 1$. We claim that if $G$ is
such that there is an exact sequence
  \begin{center}{ $G_{n-1}\rightarrow \cdots \rightarrow
G_0\rightarrow G\rightarrow 0$} \end{center}                                     with $G_{n-1}, \cdots , G_0$ all Gorenstein
injective then $G$ is Gorenstein injective. For by ([6], Proposition
1.11) $G$ is Gorenstein injective if and only if $Ext^1(L, G)=0 $ whenever
$proj.dim L<\infty$. By ([4], Corollary 4.4)
 we have that $proj.dim L<\infty$ implies $proj.dim L\leq n$.
So now using dimension shifting and these results we get that $G$ is
Gorenstein injective.
\item[(9)] If $G$ is Gorenstein injective over any Gorenstein ring $R$ and $r\in R$ is
$R$-regular, then $proj.dim R/(r)=1$. So $Ext ^1(R/(r), G)=0$ by (8). This gives
that $Hom(R,G)\stackrel{r}{\rightarrow} Hom(R,G)\rightarrow 0$ is
exact. This means that $G\stackrel{r}{\rightarrow}G $ is surjective.
So for every $x\in G$ there is a $y\in G$ with $ry=x$. Consequently we get
that $G\otimes T=0$ if $T$ has property $t(P)$ and if $r\in P$. For if $x\in G$ and $y
\in T$ and $n\geq 1$ we have that $x=r^n\overline{x}$ for some
$\overline{x}\in G$. So $x\otimes y= r^n\overline{x}\otimes y =
\overline{x}\otimes r^ny$. But if $n$ is sufficiently large we have
$r^ny=0$. Hence $x\otimes y=0$.\\
Now if $P\in Spec(R)$ and if $ht(P)\geq 1$ then since $R$ is Gorenstein (and
so Cohen-Macaulay) there is an $R$-regular element $r\in P$. Consequently we get
that $G\otimes T=0$ whenever $G$ is Gorenstein injective and when
$T$ has property $t(P)$ with $ht(P)\geq 1$.
\end{itemize}

{\bf 2. Torsion products of injective and Gorenstein injective modules}

\vspace{0.1cm}

In this section $R$ will be a Gorenstein ring of finite Krull dimension
$n$. We let $X=Spec(R)$. When we refer to $(1),(2),\cdots,(9)$ we mean the
corresponding result in the preceding section.\\

\vspace{.1cm}

{\bf Lemma 2.1.} {\it If $P\in X$ and $ht(P)\geq 1$ then for any Gorenstein
injective module $G$ we have $E(R/P)\otimes G=0$}

\vspace{0.1cm}

{\it Proof.}  By (5) we know that $E(R/P)$ has property $t(P)$. So this result
is a special case of (9).
\hfill{$\square$}
\vspace{0.1cm}

{\bf Proposition 2.2.} {\it If $G$ is Gorenstein injective and $P\in X$
then $Tor_k(E(R/P), G)=0$ if $ht(P)\neq k$. }\\

\vspace{0.1cm}

{\it Proof.} By (3) we know that $flat.dimE(R/P)=ht(P)$ so $Tor_k(E(R/P),-)=0$
if $k>ht(P)$. So we only need prove that $Tor_k(E(R/P),G)=0$ when
$G$ is Gorenstein injective and $k<ht(P)$. We prove this by induction on
$k$. If $k=0$, then $Tor_k(E(R/P), G)=E(R/P)\otimes G=0 $ if $ht(P)\geq 1$
 and $G$ is Gorenstein injective by Lemma 2.1.\\
So now we make an induction hypothesis and let $ht(P)>k$ and let $G$
be Gorenstein injective. We have an exact sequence $0\rightarrow
H\rightarrow E\rightarrow G\rightarrow 0$ with $E$ injective and $H$
Gorenstein injective. We have the exact sequence $Tor_k(E(R/P), E)
\rightarrow Tor_k(E(R/P),G)\rightarrow Tor_{k-1} (E(R/P), H)$. By the
induction hypothesis and the fact that $ht(P)>k>k-1$ we have that
$Tor_{k-1} (E(R/P), H)=0$. But $Tor_k (E(R/P), E)= 0$ by (4) and
so $Tor_k(E(R/P), G)=0$.
\hfill{$\square$}

\vspace{.1cm}

{\bf Corollary 2.3.} {\it If $0\rightarrow G'\rightarrow G\rightarrow
G''\rightarrow 0$ is an exact sequence of Gorenstein injective modules
 and if $E$ is an injective module, then for any $k\geq 0$ the
sequence $ 0\rightarrow Tor_k(E,G')\rightarrow Tor_k(E, G)
\rightarrow Tor_k( E, G'')\rightarrow 0$ is exact.}

\vspace{.1cm}

{\it Proof.}  By (1) $E$ is a direct sum of submodules isomorphic to
$E(R/P)$ with $P\in X$, it suffices to prove the claim when $E=E(R/P)$
for any $P$. In this case the claim follows from the considering the long
exact sequence of $Tor(E(R/P),-)$ associated with $0\rightarrow G'
\rightarrow G\rightarrow G''\rightarrow 0$ and Proposition 2.2.
\hfill{$\square$}

\vspace{.1cm}

{\bf Proposition 2.4.} {\it If $G$ is Gorenstein injective and $E$ is injective
then for any $k\geq 0$ the module $Tor_k(E,G)$ is a Gorenstein injective module.}

\vspace{.1cm}

{\it Proof.}  By (8) we have an exact sequence  $\cdots \rightarrow E_2
\rightarrow E_1 \rightarrow E_0 \rightarrow G\rightarrow 0$ with all the
$E_i$ injective modules where the kernels of $E_0\rightarrow G$, $E_1\rightarrow
E_0, \cdots $ are Gorenstein injective. So we can split the exact sequence
into short exact sequences $0\rightarrow G_1\rightarrow E_0 \rightarrow G \rightarrow 0$, $
0\rightarrow G_2\rightarrow E_1\rightarrow G_1\rightarrow 0,  \cdots $ with
each $G_k$ and $G$ Gorenstein injective. We then apply  Corollary 2.3 and
splice the resulting short exact sequences together to get the exact sequence
 $\cdots \rightarrow Tor_k(E,E_1)\rightarrow Tor_k(E, E_0)\rightarrow Tor_k (E, G)
\rightarrow 0$. Since each $Tor_k(E, E_n)$ is injective we get that
$Tor_k(E,G)$ is Gorenstein injective by (8).
\hfill{$\square$}

\vspace{.2cm}

{\bf 3. Filtrations of Gorenstein injective modules}

\vspace{.1cm}

We again let  $R$ be a Gorenstein ring of finite Krull dimension $n$ and
let $X=Spec(R)$ and let $X_k\subset X$ for $k\geq 0$ consist of the
$P\in X$ such that $ht(P)=k$. In this section we will also appeal to
the results $(1),(2),\cdots,(9)$ of the first section.\\
\vspace{.1cm}

The main contribution of this paper is the  following result.\\

\vspace{.1cm}

{\bf Theorem 3.1.} {\it If $G$ is a Gorenstein injective module then $G$
has a filtration  $0=G_{n+1}\subset G_n \subset \cdots \subset G_2
\subset G_1 \subset G_0=G$ where each $G_k/G_{k+1}$ ($0\leq k \leq n$) is
Gorenstein injective and has a direct sum decomposition indexed by the
$P\in X_k$ such that the summand, say $S$, corresponding to $P$ has
the property  $t(P)$ (see (5)).
Furthermore such filtrations and direct sum decompositions are unique and
functorial in $G$.}

\vspace{.1cm}

{\it Proof.} We first comment that ``functorial in $G$" means that if $H$
is another Gorenstein injective module with  such a filtration $0=H_{n+1}
\subset H_n \subset \cdots \subset H_1 \subset H_0=H$ where $T$ is the summand
of $H_k/H_{k+1}$ corresponding to $P\in X_k$ and if $f:G\rightarrow H$ is
linear then $f(G_k)\subset H_k$ for each $k$ and the induced map $G_k/G_{k+1}
\rightarrow H_k/H_{k+1}$ maps $S$ (as in the theorem) into $T$. \\
Now let $0\rightarrow R\rightarrow E^0(R)\rightarrow \cdots \rightarrow
E^n(R)\rightarrow 0$ be the minimal injective resolution of $R$ and let
$\cdots \rightarrow P_2\rightarrow P_1\rightarrow P_0\rightarrow G
\rightarrow 0$ be a projective resolution of $G$. We form the double
complex\\
\vspace{.2cm}

$
\begin{CD}
  @. 0 @. @. 0 @. \\
@. @AAA @. @AAA @. \\
0 @>>> E^0(R)\otimes P_0 @>>> \cdots \cdots \cdots @>>> E^n(R)\otimes P_0  @>>> 0\\
@. @AAA @. @AAA @.\\
0 @>>> E^0(R)\otimes P_1 @>>> \cdots \cdots \cdots @>>> E^n(R)\otimes P_1 @>>> 0\\
@. @AAA @. @AAA @.\\
 @. \vdots @. @. \vdots @.
\end{CD}
$
\vspace{.1cm}

We now use a simple spectral sequence argument. First we note that this double complex
can be regarded as a third quadrant double complex (using a shift in indices). So this will guarantee convergence of our spectral sequences. For the $E^1$ term of our
first spectral sequence we compute homology of this double complex using
the horizontal arrows.  Since each $P_n$ is projective, and so flat, we now
get the transpose of the diagram
\vspace{.2cm}

\begin{center}
{$ \cdots \rightarrow R\otimes P_1 \rightarrow R\otimes P_0 \rightarrow 0$}
\end{center}
\vspace{.2cm}
where all the missing terms are $0$. But now when we compute homology we
just get $G$ (in the $(0,0)$ position).\\
 We now first use the vertical arrows to compute
homology. The terms we get will all be of the form $Tor_i(E^j(R), G)$. By Proposition
 2.2 and (3) these are 0 unless
$i=j$. So we get a diagonal double complex. Hence the horizontal differentials
will be 0  and when we compute homology again we get $\oplus_{k=0}^n
Tor_k (E^k(R), G)$.
 This means that $G$ has a filtration $0=G_{n+1}\subset G_n
\subset \cdots \subset G_1\subset G_0=G $ with $G_k/G_{k+1} \cong
Tor_k (E^k(R), G)$ for $0\leq k \leq n$. By Proposition 2.4 we know
that each of these terms is Gorenstein injective.\\
 By (3) $E^k(R)=\oplus_{P\in X_k} E(R/P)$ and so we have that $$Tor_k (E^k(R), G)=
\oplus_{P\in X_k} Tor_k (E(R/P), G).$$ Since
each $E(R/P) $  has property $t(P)$ by (5) so does $Tor_k(E(R/P),G)$.\\
The uniqueness and functoriality will now follow from (7), i.e.
 if $P, Q$ are prime ideals of $R$ with $P\not\subset Q$ then
$Hom(S, T)=0$ whenever $S$ and $T$ have properties $t(P)$ and $t(Q)$
 respectively.\\
We now indicate how this observation gives us the functoriality and
uniqueness. Let $0\subset G_n \subset \cdots \subset G_1 \subset G$ and
$0 \subset H_n \subset \cdots \subset H_1 \subset H$ be  filtrations
of the Gorenstein injective modules $G$ and $H$ satisfying the  conclusion
of the theorem. Let $S\subset G_n $ be the summand of $G_n$ corresponding
to a given maximal ideal $P$ of $R$. Assume $n\geq 1$. Then we use the observation
that $Hom(S, U)=0$ if $U\subset H/H_1$ is the summand corresponding to
some $Q\in X_0$. Since this holds for all such $U$ we get that
$S\hookrightarrow G\rightarrow H/H_1$ is 0. So $f(S)\subset H_1$. Since this is true for all
the summands $S$ of $G_n$  we get that $f(G_n)\subset H_1$. But then
we use this argument to get $f(G_n)\subset H_2, \cdots$ and finally that
$f(G_n)\subset H_n$.\\
Repeating the argument but applied to $G/G_n\rightarrow H/H_n$ with the induced
 filtration, we get
that $f(G_{n-1})\subset H_{n-1}$ and then by the induction hypothesis that $f(G_k)\subset
H_k$ for $0\leq k\leq n$.\\
Now if $P\in X_k$ and if $S$ and $T$ are the summands of $G_k/ G_{k+1}
$ and $H_k/H_{k+1}$ corresponding to $P$ respectively then the same type argument gives
that $G_k/G_{k+1}\rightarrow H_k/H_{k+1} $ maps $S$ into $T$.\\
The uniqueness of the filtrations and direct sum decompositions can be
argued by assuming $G=H$ (with possibly different filtrations and direct
sum decompositions) and letting $f=1_G$. So the above would give
$G_k\subset H_k$. Then similarly we get $H_k\subset G_k$ and so $G_k
=H_k$ for all $k$. Likewise we get the uniqueness of the direct sum
decompositions.
\hfill{$\square$}
\vspace{.2cm}

{\bf Remark 3.2}. We would like to thank the referee for his/her help in
writing this paper.
 The referee has pointed out that the $G_k$ of Theorem 3.1
can be described by the formulas $G_k/G_{k+1}=\oplus_{P\in X_k} \Gamma_{P} (G
/G_{k+1})$ for $k=0,\cdots, n$ where for a module $M$ we have $\Gamma_P(M)$ consists of all
$x\in M$ such that $P^nx=0$ for some $n\geq 1$.
The referee also suggested that Theorem  3.1 might hold when we only assume
the ring $R$ is Cohen-Macaulay admitting a canonical module. We do not
know if this is the case.\\
\vspace{.2cm}

{\bf 4. Tensor Products of Gorenstein Injective Modules}

\vspace{.2in}

We let $R$ be a Gorenstein ring of finite Krull dimension $n$. We want to
show that over such an $R$ all tensor products of Gorenstein injective modules
 are Gorenstein injective. If $G$ (or $H$) is a Gorenstein injective
module and $0\leq k \leq n+1$ then $G_k$ (or $H_k$) will denote the {\it k}-th submodule of $G$ (or $H$)
that is part of the filtration provided by Theorem 3.1.
\vspace{.2in}

{\bf Theorem 4.1.} {\it If $G$ and $H$ are Gorenstein injective modules
then so is $G\otimes H$}.
\vspace{.2cm}

{\it Proof.} If $S$ and $T$ are Gorenstein injective $R$-modules having
properties $t(P)$ and $t(Q)$ respectively then $S\otimes T=0 $ if
$P\neq Q $ (by (5)) and if $P=Q$ and $ht(P)\geq 1$  (by (9)). We use
this to argue that $G\otimes H=G/G_1 \otimes H/H_1$. This claim
is trivial if $n=0$ since then $G_1=H_1=0$. So suppose $n\geq 1$.
Then using the above and Theorem 3.1 we see that $G_n\otimes H_k/H_{k+1}
=0 $ for $k=0,\cdots, n$. Hence $G_n\otimes H=0$. Then tensoring
the exact sequence $0\rightarrow G_n\rightarrow G\rightarrow G/G_n\rightarrow 0$
with $H$ we get that $G\otimes H=G/G_n \otimes H$.\\
 If $n\geq 2 $ (i.e. $n-1\geq 1)$ then the same argument gives that
$G_{n-1}/G_n \otimes H=0$ and then that $G\otimes H=G/G_{n-1} \otimes H$.
Continuing we get that $G\otimes H= G/G_1 \otimes H$. But then the
same type argument gives that $G/G_1\otimes H=G/G_1 \otimes H/H_1$
and so that $G\otimes H=G/G_1 \otimes H/H_1$.\\
Now by Theorem 3.1 and (5) we see that $G\otimes H=G/G_1 \otimes H/H_1$
will be a direct sum of modules of the form $S\otimes T$ where $S$ and
$T$ both have property $t(P)$ for a minimal prime ideal $P$ of $R$.
But such an $S$ and $T$ are naturally modules over $R_P$ and hence
$S\otimes T$ is an $R_P$-module. Then by (8) $S\otimes T$ is a Gorenstein
injective $R$-module. So finally noting that the class of Gorenstein
injective modules is closed under arbitrary direct sums (by (8))
we get that $G\otimes H$ is a Gorenstein injective
$R$-module. \hfill{$\square$}\\

\vspace{.2cm}

{\bf Remark 4.2}. With the same hypothesis as in the Theorem 4.1, we
do not know if each $Tor_k (G,H)$ is also Gorenstein injective when
$k>0$.\\

\vspace{.2cm}

{\bf Acknowledgements.} This paper was written while the second author
was visiting the University of Kentucky from June to August, 2009. The
second author was partially supported by the Specialized Research Fund
for the Doctoral Program of Higher Education (Grant No. 20060284002) and
NSFC (Grant No. 10771095) and NSF of Jaingsu Province of China.

\end{document}